\documentclass[11pt,reqno]{amsart}

\pdfoutput = 1

\usepackage[english]{babel}
\usepackage{amsmath,amssymb,amsthm}
\usepackage{amsfonts}
\usepackage{amsxtra}

\usepackage{array,dcolumn}
\usepackage{graphicx}
\usepackage[hyperfootnotes=true,
        pdffitwindow=true,
        plainpages=false,
        pdfpagelabels=true,
        pdfpagemode=UseOutlines,
        pdfpagelayout=SinglePage,
                hyperindex,]{hyperref}

\usepackage[T1]{fontenc}
\usepackage[utf8]{inputenc}
\usepackage[activate={true,nocompatibility},spacing,kerning]{microtype}
\usepackage{bm}
\usepackage{bbm}
\usepackage[sc,osf]{mathpazo}
\microtypecontext{spacing=nonfrench}

 \calclayout

\newcommand{\mathsym}[1]{{}}
\newcommand{\unicode}[1]{{}}

\theoremstyle{plain}

\theoremstyle{definition}

\theoremstyle{remark}

\renewcommand{\leq}{\leqslant}
\renewcommand{\geq}{\geqslant}

\numberwithin{equation}{section}

\begin{document}

\title{Comment on "Sum of squares of uniform random variables" by I.~Weissman}
\author{Peter J. Forrester}
\address{ARC Centre of Excellence for Mathematical and Statistical Frontiers,
School of Mathematics and Statistics, The University of Melbourne, Victoria 3010, Australia.} 
\email{pjforr@unimelb.edu.au}
\date{}
%

\begin{abstract}
\noindent
The recent paper by I.~Weissman, "Sum of squares of uniform random variables",
[Statist. Probab. Lett. \textbf{129} (2017), 147--154] is compared to earlier work of
B.~Tibken and D.~Constales relating to the area of the intersection of a centred ball and
cube in $\mathbb R^n$, published in the  Problems and Solutions section of SIAM Review in
1997. Some recent applications of explicit formulas 
for the corresponding probabilities from these references, to problems in lattice reduction, and to
the study of Lyapunov exponents of products of random matrices, are noted.
\end{abstract}

\maketitle

\section{Introduction}\label{s1}
A recent paper in this journal by Weissman \cite{We17} addresses the question of the distribution of the sum of 
the squares of $n$ independent uniform random variables in $[0,1]$. It is stated in the Introduction
of  \cite{We17} that there is nothing
in the published literature dealing with this problem. While this may be true in literal terms, the solution
to an equivalent question within geometrical probability has been available in the literature since 1997,
in the form of two answers to a particular problem posed in the section "Problems and Solutions" of
SIAM Review \cite{RR97}.

The purpose of this note is to explicitly connect \cite{RR97} with the topic of \cite{We17}, to isolate from
 \cite{RR97}  formulas which supplement those obtained in  \cite{We17}, and to indicate some recent
 developments in the form of both applications and new theory.

\section{Comparison of results}\label{s1a} 
Let $U_i$, $1 \le i \le N$ be independent, uniformly distributed random variables on $[0,1]$, and
set 
\begin{equation}\label{15.1}
S_N = \sum_{i=1}^N U_i^2. 
\end{equation}
With $\chi_J = 1$ for $J$ true, 
$\chi_J=0$ otherwise, one sees immediately that \cite[Eq.~(7)]{We17}
\begin{align}\label{FV}
F_N(s) & := {\rm Pr} \, (S_N \le s) = \int_0^1 dx_1 \cdots \int_0^1 dx_N \,
\chi_{\sum_{j=1}^N x_j^2 \le s} \nonumber \\
& = {\rm Vol} \, ( C_N([0,1]) \cap B_N(\sqrt{s}) ),
\end{align}
where $C_N([a,b])$ denotes the cube in $\mathbb R^N$, $[a,b]^N$, and $B_N(r)$ denotes the 
ball in $\mathbb R^N$, centred at the origin and of radius $r$. By symmetry and scaling this can
be written
\begin{align}\label{FU}
F_N(s) & = 2^{-N} {\rm Vol} \, ( C_N([-1,1]) \cap B_N(\sqrt{s}) ) \nonumber \\
& = 2^{-N} s^{N/2} {\rm Vol} \, ( C_N([-1/\sqrt{s},1/\sqrt{s}]) \cap B_N(1) ).
\end{align}

The Problems and Solutions section of volume 39 of SIAM Review, 
published in 1997 \cite{RR97}, contained solutions by
B.~Tibken and D.~Constales of Problem 96-19 as posed by Liqun Xu:
"What is the volume of the intersection of a cube and a ball in $N$-space? Assume that the cube
and the ball have a common centroid and that neither of the two bodies is large enough to contain
the other completely? This problem arises in the study of certain probability distributions."
One sees that (\ref{FU}) is precisely the geometric quantity as sought by Xu.

The explicit form for $N=2$,
\begin{equation}\label{2.a}
F_N(s) = \left \{
\begin{array}{ll} {1 \over 4} \pi s, & 0 \le s \le 1 \\
\sqrt{s - 1} + (\pi/4 -  {\rm arccos} \, (1/\sqrt{s})) s, & 1 \le s \le 2 \\
1, & 2 \le s. \end{array} \right.
\end{equation}
was derived in both the solutions of Tibken and Constales (for the latter use needs to be made of the
second equality in (\ref{FU})). This is \cite[Eq.~(10)]{We17}. In relation to $N=3$, 
let
$$
h(a) = 8 a^2 \sqrt{1 - 2 a^2} +
2(3a - a^3) \Big ( 4 {\rm arcsin} \, \Big ( {a \over \sqrt{1 - a^2}} \Big ) - \pi \Big ) -
8 {\rm arcsin} \, \Big ( {a^2 \over 1 - a^2} \Big ) + {4 \pi \over 3}.
$$
From the solution of Tibken, Eq.~(23), we read off that
\begin{equation}\label{2.b}
 {\rm Vol} \, ( C_N([-a,a]) \cap B_N(1) )  =
 \left \{
\begin{array}{ll} 8a^3, & 0 \le a < 1/\sqrt{3} \\
h(a), & 1/\sqrt{3} \le a < 1/\sqrt{2} \\
\pi(6a - 2a^3 - 8/3), & 1/\sqrt{2} \le a \le 1 \\
4 \pi/3, & 1 \le a. \end{array} \right.
\end{equation}
Recalling the 2nd equality in (\ref{FU}), we see this is equivalent to \cite[Eq.~(13)]{We17}.

For general $N$, using a method based on Fourier series, it was proved by Constales\footnote{
after the correction $2 \sqrt{k/N} \mapsto \sqrt{2\pi k/N}$ for all quantities relating to the power of
$N$ in the sum, as alerted to me by E.~Postlethwaite, and as is consistent with \cite[Th.~4]{AN17}}
 that
\begin{equation}\label{C1}
F_N(s) = \bigg ( {1 \over 6} + {s \over N} + {1 \over \pi} {\rm Im} \,
\sum_{k=1}^\infty \Big ( {C( \sqrt{2 \pi k/N}) - i S( \sqrt{2 \pi k/N}) \over
 \sqrt{2 \pi k/N}} \Big )^N {e^{2 \pi i k s/N} \over k}  \bigg ).
\end{equation}
Here $S(x) = \int_0^x (\sin t^2) \, {\rm d}t$ and  $C(x) = \int_0^x (\cos t^2) \, {\rm d}t$ denote the
Fresnel integrals. Let
${\mathcal L}(F_N)(s):= \int_0^\infty F_N(t) e^{-st} \, {\rm d}t$ denote the Laplace transform of $F_N$.
Constales also remarks that
$$
{\mathcal L}(F_N)(s) = {2^{-N} \pi^{N/2} ({\rm erf} \, \sqrt{s})^N \over s^{N/2+1}},
$$
and thus by the inversion formula for the Laplace transform
\begin{equation}\label{C2}
F_N(t) = 2^{-N} {\pi^{N/2} \over 2 \pi i}
\int_{c - i \infty}^{c + i \infty}
{ ({\rm erf} \, \sqrt{s})^N \over s^{N/2 + 1} } e^{st} \, {\rm d}s, \qquad c> 0.
\end{equation}
It is commented too that the central limit theorem applied to (\ref{15.1})
implies that for large $N$, $F_N'(s)$ has the form of the normal distribution
${\rm N}[N/3, 2 \sqrt{N /{45}}  ]$ cf.~the first sentence of
\cite[\S2]{We17}.

Using a method based on the Fourier integral rather than Fourier series, it was proved by
Tibken that for general $N$
\begin{equation}\label{15.2a}
F_N(s) = s^{N/2} \Big ( {1 \over 2 \pi} \Big )^{N/2}
\int_{\mathbb R^N}
{J_{N/2}((\sum_{j=1}^N y_j^2)^{1/2}) \over
(\sum_{j=1}^N y_j^2)^{N/4}}
\prod_{j=1}^N {\sin ( y_j/\sqrt{s}) \over y_j} \, {\rm d}y_1 \cdots {\rm d}y_N
\end{equation}
(this is Eq.~(13) in Tibken with use required too of (\ref{FU}) since Tibken computes
the quantity Vol$\,(C_N([-a,a]) \cap B_N(1))$).
Here $J_p(z)$ denotes the Bessel function of order $p$. Moreover, with
$$
I_a(\lambda,0) = \Big ( \pi \, {\rm erf} \Big ( {a \over 2 \sqrt{\lambda}} \Big ) \Big )^n,
\qquad
I_a(\lambda,k) = (-1)^k {\partial^k \over \partial \lambda^k} I_a(\lambda,0),
$$
and $L_n^{(k)}(x)$ denoting the classical Laguerre polynomial, it was shown in Eq.~(33)
of Tibken's work that (\ref{15.2a}) can
be reduced to
\begin{equation}\label{T1}
2^N s^N F_N(1/s^2) = \pi^{-N/2} \bigg (
{I_s(1/(2N+4),0) \over \Gamma(N/2+1)} +
\sum_{k=2}^\infty {L_k^{N/2}(N/2+1) I_s(1/(2N+4),k) \over
\Gamma(k+N/2+1) (2N+4)^k} \bigg ).
\end{equation}

\section{Applications}
\subsection{Complex lattice reduction}\label{S3.1}
Let $\mathcal B = \{ \mathbf b_1, \mathbf b_2,\dots, \mathbf b_n \}$ be a basis in
$\mathbb R^n$. Associated with $\mathcal B$ is the lattice
\begin{equation}\label{L1}
\mathcal L_{\mathcal B} = \Big \{ n_1 \mathbf b_1 + n_2 \mathbf b_2 +
\cdots  + n_N \mathbf  b_N \, \Big | \, n_i \in \mathbb Z \: (i=1,\dots,N) \Big \}.
\end{equation}
For dimensions $N \ge 2$, there are an infinite number of bases giving the same lattice.
A basic question is to determine a basis --- referred to as a reduced basis --- consisting of 
short and/ or  near orthogonal vectors. In two dimensions the Lagrange--Gauss lattice
reduction algorithm (see e.g.~\cite{Br12}) performs this task, and converges to basis vectors
$\{ \bm{\alpha}, \bm{\beta} \}$ with the property that
\begin{equation}\label{3.1}
|| \bm{\alpha} || \le || \bm{\beta} ||, \qquad
\Big | {\bm{\alpha} \cdot \bm{\beta} \over || \bm{\alpha} ||^2 } \Big | \le {1 \over 2}.
\end{equation}

It is well known, and straightforward to verify, that the inequalities (\ref{3.1}) imply that 
these basis vectors are the shortest possible.
Probabilities and volumes can be introduced into this setting by introducing
the notion of a random basis, as then the reduced basis is given by a probability distribution. 
The most natural meaning of a random basis, normalised to have a unit cell with volume unity, is to
choose the matrix $M$ of lattice vectors from ${\rm SL}_N(\mathbb R)$, endowed with a
Haar measure as identified by Siegel \cite{Si45}. 

The problem of computing (\ref{FV}) shows itself for $N=2$ in this context by considering the
lattice reduction problem for a complex analogue of (\ref{L1}) in two-dimensions.
Specifically, one requires that $\mathbf b_1$, $\mathbf b_2$ be linearly independent vectors in
$\mathbf C^2$, and chooses $n_1, n_2 \in \mathbb Z[i]$, and thus as Gaussian integers.
In a QR (Gram-Schmidt) parametrisation, the matrix of basis vectors is decomposed $M = U T$, where
$U \in {\rm SU}(2)$ and
$$
T = \begin{bmatrix} t_{11} & t_{12}^{(r)} + i t_{12}^{(i)} \\
0 & t_{22} \end{bmatrix}, \qquad t_{11} > 0, \quad t_{22} = 1/t_{11}.
$$
Here we think of $U$ as rotating the lattice so that $\mathbf b_1 = (r_{11}, 0)$,
$\mathbf b_2 = (t_{12}^{(r)} + i t_{12}^{(i)}, t_{22})$. After integrating over the variables
associated with $U$, and $t_{22}$, and introducing $t_{11} = t$, $t_{12}^{(r)} = y_1$,
$t_{12}^{(i)} = y_2$, for notational convenience, the invariant measure restricted to the
domain of the shortest reduced basis is equal to \cite[Eq.~(4.26)]{FZ17}
\begin{equation}\label{24.1}
(2\pi^2)t\chi_{\lVert\mathbf{y}\rVert^2\geq t^2-1/t^2}\chi_{|y_1|\leq t/2}\chi_{|y_2|\leq t/2}\mathrm{d}t\mathrm{d}y_1\mathrm{d}y_2.
\end{equation}
As noted in \cite{FZ17}, with $V_2(a,b)$ denoting the area of overlap between a disk of radius $a$, and 
a square of side length $b$, both centred at the origin, (\ref{24.1}) upon integration over
$y_1$ and $y_2$ can be written
\begin{equation}\label{24.1a}
(2\pi^2)\bigg(t^3\chi_{0<t<1}+\chi_{t>1}t\left(t^2-V_2\left((t^2-1/t^2)^{1/2},t/2\right)\right)\bigg)\mathrm{d}t.
\end{equation}
Use of (\ref{2.a}) and (\ref{FU}) in (\ref{24.1a}) allows for the determination of the explicit functional
form of $V_2$, which up to normalisation corresponds to the PDF of the length of the shortest lattice
vector. The latter quantity can be realised as a numerically generated histogram,
by Monte Carlo sampling from SL${}_2(\mathbb C)$
with Haar measure, and then use of the complex Lagrange--Gauss lattice reduction algorithm to determine
the shortest basis for each sample;
see \cite[\S 6]{FZ17}.

\subsection{Lyapunov exponents}
Let $\mathbf r_i \in \mathbb R^3$ $(i=1,2,3)$ be sampled uniformly from the unit sphere in
$\mathbb R^3$; this can be done for example by setting $\mathbf r_i = (x_i, y_i, z_i)/
\sqrt{x_i^2 + y_i^2 + z_i^2}$, where each of $x_i,y_i,z_i$ are independent standard
normal random variables. Specify an element of the ensemble of $2 \times 2$ random matrices
$U_2^B$ by forming $\begin{bmatrix} x_1 & y_1 \\ x_2 & y_2 \end{bmatrix}$, and thus restricting
$\mathbf r_1$, $\mathbf r_2$ to the first two components along each row.
Specify an element of the ensemble of $3 \times 3$ random matrices $U_3^S$ by
forming
$$
\begin{bmatrix} x_1 & y_1 & z_1 \\
x_2 & y_2 & z_2 \\
x_3 & y_3 & z_3
\end{bmatrix},
$$
and so placing the components of $\mathbf r_i $ along each row in order.
The ensembles  $U_2^B$ and $U_3^S$ both have the property that the distribution of
any member, $X_i$ say, is unchanged by multiplication by an appropriately sized
real orthogonal matrix on the right. As such, the Lyapunov exponent
$\mu_1$, defined as
\begin{equation}\label{D1}
\mu_1 = \lim_{m \to \infty} {1 \over m} \log
{\rm max}_{|| \mathbf v|| = 1} || X_m X_{m-1} \cdots X_1 \mathbf v ||,
\end{equation}
is given by the simple formula \cite{CN84}
$$
2 \mu_1 = \int_0^\infty (\log t) p(t) \, {\rm d}t,
$$
where $p(t)$ is the PDF for $\sum_{i=1}^p x_i^2$, ($p=2,3$), corresponding to  the sum 
of the squares of the entries in
the first column of any $X_i$.

We know from \cite[Prop.~6]{FZ18} that both the ensembles 
$U_2^B$ and $U_3^S$ have each $x_i$ uniformly distributed on $[-1,1]$
and thus $p(s) = F_N'(s)$ for $N=2$, $N=3$ respectively.
Consequently for $U_2^B$, using (\ref{2.a})
\begin{align}\label{M1}
2 \mu_1 & = {\pi \over 4} \int_0^1 \log s \, {\rm d}s +
\int_1^2 \Big ( {\rm arcsin} \Big ( {1 \over \sqrt{s}} \Big ) - {\pi \over 4} \Big ) \, \log s \, {\rm d}s \nonumber \\
& \approx -0.736056.
\end{align}
And for $U_3^S$, as implied by (\ref{2.b}) and (\ref{FU}), or more explicitly by
\cite[Eq.~(14)]{We17}, with the notation
\begin{align*}
f_{3,2}(s)  & = 3 \Big ( {\rm arcsin} \Big ( {1 \over \sqrt{s -1}} \Big ) - {\pi \over 4} 
\Big ) +
\sqrt{s} \bigg ( {\rm arctan} \, \sqrt{s-2 \over s} -  {\rm arctan} \, \sqrt{1 \over s(s-2)} \bigg ), \\
& =  3 \Big ( {\rm arcsin} \Big ( {1 \over \sqrt{s -1}} \Big ) - {\pi \over 4} 
\Big ) + \sqrt{s} \bigg ( {\pi \over 4} - {3 \over 2}  {\rm arcsin} \Big ( {1 \over {s -1}} \Big ) \bigg )
\end{align*}
we have
\begin{align}\label{M2}
2 \mu_1 & = {\pi \over 4} \int_0^1 s^{1/2} \log s \, {\rm d}s +
{\pi \over 4} \int_1^2 (3 - 2 s^{1/2}) \log s \, {\rm d}s +
\int_2^3 f_{3,2}(s) \log s \, {\rm d}s
 \nonumber \\
& \approx   -0.187705.
\end{align}
Numerical methods similar to those used in \cite[\S 4]{Fo12}, which are
based on the definition  (\ref{D1}), can be used to obtain
Monte Carlo estimates of (\ref{M1}) and (\ref{M2}), and so realising the numerical
values to several decimal places.

\subsection{Random sampling --- lattice enumeration with discrete pruning}
In subsection \ref{S3.1}, attention was drawn to the problem of determining a
reduced basis for a general lattice (\ref{L1}). A more modest aim, but one with wide
applicability --- for example to lattice based cryptography \cite{Ng09} --- is to find a
nonzero lattice vector of the smallest norm. This is the so-called shortest vector problem.
In the enumeration approach to the shortest vector problem
\cite{FP85}, the basic idea is to perform a depth-first search on a tree whose leaves correspond
to lattice points, and internal nodes to coefficients of the integer combination
 practical improvement is to restrict the exhaustive search to a subset of possible solutions,
 by pruning subtrees of the tree for which the "probability" of finding the shortest
 vector is small \cite{SE94}. In \cite{AN17} a different pruning set to that in \cite{SE94} is
 introduced. Analysis of the algorithm requires computing
 $$
 {\rm Vol} \, \Big ( \prod_{l=1}^N [ \alpha_l, \beta_l] \cap B_N(1) \Big ), \qquad
  \prod_{l=1}^N [ \alpha_l, \beta_l] :=
  \{ (x_1,\dots, x_N): \, x_l \in [\alpha_l, \beta_l] \: (l=1,\dots,N) \}
  $$
  (cf.~(\ref{FU})). Aono and Nguyen \cite{AN17} give generalisations of each of 
  Constales formulas (\ref{C1}), (\ref{C2}), and Tibken's formula (\ref{T1}). The
  simplest is the generalisation of (\ref{C2}), which reads
  \begin{equation}\label{AN1}
  {\rm Vol} \, \, \Big ( \prod_{l=1}^N [ \alpha_l, \beta_l] \cap B_N(1) \Big ) =
  {\pi^{N/2} \over 2 \pi i}
  \int_{c - i \infty}^{c + i \infty}
  {e^s \over s^{N/2 + 1}} 
  \prod_{j=1}^N {({\rm erf} (\beta_j \sqrt{s}) - {\rm erf} (\alpha_j \sqrt{s})) \over
  2 (\beta_j - \alpha_j)} \, {\rm d}s,
  \end{equation}
  for suitable $c$. The efficient numerical computation of (\ref{AN1}) is discussed in
  detail in \cite{AN17}.

  \section*{Acknowledgements}
This research project is part of the program of study supported by the 
ARC Centre of Excellence for Mathematical \& Statistical Frontiers,
and the Australian Research Council Discovery Project grant DP170102028.
The referee of the published form of this work [Statist. Probab. Lett. \textbf{129}
(2018), 147--154] is to be thanked for picking up a number of errors relative to the original
arXiv posting. However, this report was unfortunately missed at the time, and not acted
on until now. Most recently (11${}^{\rm th}$ August, 2022) Eamonn Postlethwaite
alerted me to the error in the form of (\ref{C1}) as reported in Constales
article from 1997, which was copied into my original arXiv posting, and is now
fixed.

\providecommand{\bysame}{\leavevmode\hbox to3em{\hrulefill}\thinspace}
\providecommand{\MR}{\relax\ifhmode\unskip\space\fi MR }
\providecommand{\MRhref}[2]{%
  \href{http://www.ams.org/mathscinet-getitem?mr=#1}{#2}
}
\providecommand{\href}[2]{#2}

\end{document}